\documentclass{article}
\usepackage{ijcai22}

\usepackage{times}
\usepackage{soul}
\usepackage{url}
\usepackage[hidelinks]{hyperref}
\usepackage[utf8]{inputenc}
\usepackage[small]{caption}
\usepackage{graphicx}
\usepackage{amsmath}
\usepackage{amsthm}
\usepackage{booktabs}
\usepackage{algorithm}
\usepackage{algorithmic}
\title{Accelerated Optimization on Riemannian Manifolds
	\\ via Projected Variational Integrators}

\author{
    Valentin Duruisseaux and Melvin Leok
}

\graphicspath{{figures/}}

\usepackage{amssymb}
\usepackage{enumitem}
\usepackage{mathtools}
\usepackage[all]{xy}

\usepackage{MnSymbol}

\usepackage{dutchcal}

\DeclareMathOperator*{\argmin}{argmin}

\newtheorem{definition}{Definition}

\begin{document}

\maketitle

\begin{abstract}
	A variational formulation of accelerated optimization on normed spaces was recently introduced by considering a specific family of time-dependent Bregman Lagrangian and Hamiltonian systems whose corresponding trajectories converge to the minimizer of the given convex function at an arbitrary accelerated rate of $\mathcal{O}(1/t^p)$. This framework has been exploited using time-adaptive geometric integrators to design efficient explicit algorithms for symplectic accelerated optimization. It was observed that geometric discretizations were substantially less prone to stability issues, and were therefore more robust, reliable, and computationally efficient. More recently, this variational framework has been extended to the Riemannian manifold setting by considering a more general family of time-dependent Bregman Lagrangian and Hamiltonian systems on Riemannian manifolds. It is thus natural to develop time-adaptive Hamiltonian variational integrators for accelerated optimization on Riemannian manifolds. In the past, Hamiltonian variational integrators have been constructed with holonomic constraints, but the resulting algorithms were implicit in nature, which significantly increased their cost per iteration. In this paper, we will test the performance of explicit methods based on Hamiltonian variational integrators combined with projections that constrain the numerical solution to remain on the constraint manifold.
\end{abstract}

\section{Introduction}

Many data analysis and machine learning algorithms are designed around the minimization of a loss function or the maximization of a likelihood function. Due to the ever-growing size of data sets, there has been a lot of focus on first-order optimization algorithms because of their low cost per iteration. Nesterov's accelerated gradient method~\cite{Nes83} was shown to converge in $\mathcal{O}(1/k^2)$ to the minimum of the convex objective function $f$ at hand, improving on the $\mathcal{O}(1/k)$ convergence rate exhibited by standard gradient descent methods. This phenomenon in which an algorithm displays this improved rate of convergence is referred to as acceleration. Nesterov's algorithm was shown in~\cite{SuBoCa16} to limit to a second-order ODE as the time-step goes to 0, and $f(x(t))$ converges to its optimal value at a rate of $\mathcal{O}(1/t^2)$ along the trajectories of this ODE. It was later shown that in continuous time, an arbitrary convergence rate $\mathcal{O}(1/t^p)$ can be achieved in normed spaces~\cite{WiWiJo16} and on Riemannian manifolds~\cite{Duruisseaux2021Riemannian}, by considering flow maps generated by a family of time-dependent Bregman Lagrangian and Hamiltonian systems which is closed under time-rescaling. This variational framework and the time-rescaling property of this family was then exploited in \cite{duruisseaux2020adaptive} using time-adaptive geometric integrators to design efficient explicit algorithms for accelerated optimization on normed vector spaces. It was observed that a careful use of adaptivity and symplecticity could result in a significant gain in computational efficiency. More generally, when applied to Hamiltonian systems, symplectic integrators yield discrete approximations of the flow that preserve the symplectic two-form~\cite{HaLuWa2006}, and results in the preservation of many qualitative aspects of the underlying dynamical system and, in particular, exhibit excellent long-time near-energy preservation~\cite{Re1999,Benettin1994}. Variational integrators provide a systematic method for constructing symplectic integrators of arbitrarily high-order based on the discretization of Hamilton's principle~\cite{MaWe2001}. 

Recently, there has been some effort to derive accelerated optimization algorithms in the Riemannian manifold setting~\cite{Duruisseaux2021Constrained,Duruisseaux2021Riemannian,Alimisis2020-1,Sra2016,Sra2018,Sra2020,Liu2017}. The Whitney Embedding Theorems \cite{Whitney1944_2,Whitney1944_1} state that any smooth manifold of dimension $n \geq 2$ can be embedded in $\mathbb{R}^{2n}$ and immersed in $\mathbb{R}^{2n-1}$, and is thus diffeomorphic to a submanifold of $\mathbb{R}^{2n}$. Furthermore, the Nash Embedding Theorem \cite{Nash1956} states that any Riemannian manifold can be globally isometrically embedded into some Euclidean space. As a consequence, the study of Riemannian manifolds can in principle be reduced to the study of submanifolds of Euclidean spaces. Altogether, this motivates the introduction of time-adaptive variational integrators on Riemannian manifolds that exploit the structure of the embedding Euclidean space. The time-adaptive approach relying on a Poincar\'e transformation from~\cite{duruisseaux2020adaptive} was extended to the Riemannian manifold setting in~\cite{Duruisseaux2021Riemannian}, and~\cite{Duruisseaux2021Constrained} studied how holonomic constraints can be incorporated into variational integrators to constrain the numerical solution to the Riemannian manifold. Although these integrators were carefully justified geometrically as coming from discrete action principles, they were implicit in nature, which significantly increases their cost per iteration as the dimension of the problem becomes large. 

In this paper, we present new algorithms based on explicit variational integrators in the embedding space where the manifold constraints are enforced via projections. The resulting explicit algorithms are then used to numerically solve generalized eigenvalue and Procrustes problems on the unit sphere and Stiefel manifold. We believe that these algorithms are the most efficient methods to date which exploit the variational framework from \cite{Duruisseaux2021Riemannian}.

\section{Preliminaries}

\subsection{Variational Integration} \label{SectionBasicTheory}

Variational integrators are derived by discretizing Hamilton's principle, instead of discretizing Hamilton's equations. As a result, variational integrators are symplectic, preserve many invariants and momentum maps, and have excellent long-time near-energy preservation~\cite{MaWe2001}. Traditionally, variational integrators have been designed based on the Type I generating function known as the discrete Lagrangian, but more recently, variational integrators have been extended to the framework of Type II/III generating functions, commonly referred to as discrete Hamiltonians~\cite{LaWe2006,LeZh2011,ScLe2017}. The boundary-value formulation of the exact Type II generating function of the time-$h$ flow of Hamilton's equations is given by the exact discrete right Hamiltonian,
\begin{equation*}
	H_d^{+,E}(q_0,p_h) =  p_h q_h - \int_0^h \left[ p(t) \dot{q}(t) -H(q(t), p(t)) \right] dt, \label{exact_Hd}
\end{equation*}
where $(q,p)$ satisfies Hamilton's equations with boundary conditions $q(0)=q_0$ and $p(h)=p_h$. A Type II Hamiltonian variational integrator is constructed by using an approximate discrete Hamiltonian $H_d^+$, and applying the discrete right Hamilton's equations,
\begin{equation}\label{Discrete Right Eq}
	p_0=D_1H_d^+(q_0,p_1), \qquad q_1=D_2H_d^+(q_0,p_1),
\end{equation}
which implicitly defines the discrete right Hamiltonian map $\tilde{F}_{H_d^+}:(q_0,p_0) \mapsto (q_1,p_1)$. Theorem 2.2 from~\cite{ScLe2017} states that if a discrete right Hamiltonian $H^+_d$ approximates the exact discrete right Hamiltonian $H_d^{+,E}$ to order $r$, then the discrete right Hamiltonian map $\tilde{F}_{H^+_d}$, viewed as a 1-step method, is order $r$ accurate.

\subsection{Riemannian Geometry}

We first introduce a few main notions from Riemannian geometry. See \cite{Absil2008,Boumal2020,Duruisseaux2021Riemannian,Lee2019} for more details on Riemannian manifolds and optimization on manifolds.


\begin{definition} \label{def: fiber metric}
	Let  $\mathcal{Q}$ be a Riemannian manifold with Riemannian metric  $g(\cdot,\cdot) = \langle \cdot , \cdot \rangle$. We define the \textbf{musical isomorphism} $g^{\flat}:T\mathcal{Q} \rightarrow T^*\mathcal{Q}$ via $ g^{\flat}(u)(v) = g_q(u,v) $ for all $q\in \mathcal{Q} $ and $  u,v\in T_q\mathcal{Q}$, and its \textbf{inverse musical isomorphism} $g^{\sharp}:T^*\mathcal{Q} \rightarrow T\mathcal{Q}$. The Riemannian metric $g(\cdot,\cdot) = \langle \cdot , \cdot \rangle$ induces a \textbf{fiber metric} $g^*(\cdot ,\cdot) =   \llangle \cdot , \cdot \rrangle $ on $T^* \mathcal{Q}$ via
	\[  \llangle u , v \rrangle = \langle g^{\sharp}(u), g^{\sharp}(v) \rangle  \quad \forall u,v \in T^* \mathcal{Q}. \]
\end{definition}


\begin{definition}
	Denoting the exterior derivative of $f$ by $df$, the \textbf{Riemannian gradient} $\emph{grad}f(q) \in T_q  \mathcal{Q}$ at $q\in \mathcal{Q}$ of a smooth function $f:\mathcal{Q} \rightarrow \mathbb{R}$ is the tangent vector at $q$ such that 
	\[ \langle  \emph{grad}f(q) , u \rangle = df(q) u   \qquad \forall u\in T_q \mathcal{Q}. \]
\end{definition}

\begin{definition}
	A \textbf{geodesic} in a Riemannian manifold $\mathcal{Q}$ is a parametrized curve $\gamma : [0,1] \rightarrow \mathcal{Q}$ which is of minimal local length, and is a generalization of the notion of a straight line from Euclidean spaces to Riemannian manifolds. 
	
\end{definition}

\begin{definition}
	The \textbf{Riemannian Exponential} $\emph{Exp}_q:T_q \mathcal{Q} \rightarrow \mathcal{Q}$ at $q\in \mathcal{Q}$ is defined via
$ \emph{Exp}_q(v) = \gamma_v(1), $
	where $\gamma_v$ is the unique geodesic in $\mathcal{Q}$ such that $\gamma_v(0) = q$ and $\gamma_v'(0) = v$, for any $v\in T_q \mathcal{Q} $.  $\emph{Exp}_q$ is a diffeomorphism in some neighborhood $ U \subset T_q\mathcal{Q}$ containing 0, so its inverse, the \textbf{Riemannian Logarithm} $\emph{Log}_p : \emph{Exp}_q(U) \rightarrow T_q \mathcal{Q}$, is well-defined.
\end{definition}

\begin{definition}
	A \textbf{retraction} on a manifold $\mathcal{Q}$ is a smooth mapping $\mathcal{R}: T\mathcal{Q} \rightarrow \mathcal{Q}$, such that for any $q \in \mathcal{Q}$, the restriction $\mathcal{R}_q : T_q\mathcal{Q} \rightarrow \mathcal{Q} $ of $\mathcal{R}$ to $T_q\mathcal{Q} $ satisfies
	\begin{itemize}
		\item $\mathcal{R}_q(0_q) = q$, where $0_q$ denotes the zero element of $T_q\mathcal{Q} $, 
		\item $T_{0_q}\mathcal{R}_q = \mathbb{I}_{T_q\mathcal{Q} }$ with the canonical identification  $T_{0_q}T_{q}\mathcal{Q} \simeq T_{q}\mathcal{Q}$, where $T_{0_q}\mathcal{R}_q$ is the tangent map of $\mathcal{R}$ at $0_q \in T_{q}\mathcal{Q}$ and $\mathbb{I}_{T_q\mathcal{Q} }$ is the identity map on $T_{q}\mathcal{Q}$.
	\end{itemize}
	The Riemannian Exponential map is a natural example of a retraction on a Riemannian manifold. 
\end{definition}

\begin{definition}
	A subset $A$ of a Riemannian manifold $\mathcal{Q}$ is called \textbf{geodesically uniquely convex} if every two points of $A$ are connected by a unique geodesic in $A$. A function $f:\mathcal{Q} \rightarrow \mathbb{R}$ is called \textbf{geodesically convex} if for any two points $q,\tilde{q} \in \mathcal{Q}$ and a geodesic $\gamma$ connecting them, 
	\[ f(\gamma(t)) \leq (1-t) f(q) +t f(\tilde{q})  \qquad \forall t\in [0,1].\] Note that if $f$ is a smooth geodesically convex function on a geodesically uniquely convex subset $A$,
	\[ f(q) - f(\tilde{q}) \geq \langle  \emph{grad}f(\tilde{q}) , \emph{Log}_{\tilde{q}}(q) \rangle   \qquad \forall q,\tilde{q} \in A.\]
	A function $f:A\rightarrow \mathbb{R}$ is called \textbf{geodesically $\lambda$-weakly-quasi-convex} with respect to $q \in \mathcal{Q}$ for some $\lambda \in (0,1]$ if
	\[ \lambda \left(f(q) - f(\tilde{q})\right) \geq \langle  \emph{grad}f(\tilde{q}) , \emph{Log}_{\tilde{q}}(q) \rangle   \qquad \forall \tilde{q} \in A.\]
	Since a geodesically convex function is $\lambda$-weakly-quasi-convex with $\lambda=1$, the algorithms introduced in this paper can also be used in the geodesically convex case. Note that a local minimum of a geodesically convex or $\lambda$-weakly-quasi-convex function is also a global minimum. 
\end{definition}

\begin{definition}
	Given a Riemannian manifold $\mathcal{Q}$ with sectional curvature bounded below by $K_{\min}$, and an upper bound $D$ for the diameter of the considered domain, define
	\begin{equation}\label{eq: zeta}
		\zeta = 
		\begin{cases}
			\sqrt{-K_{\min}} D \coth{ (\sqrt{-K_{\min}} D) }  & \quad \text{if }  K_{\min} < 0 \\ 1   & \quad \text{if }  K_{\min} \geq 0 
		\end{cases} .
	\end{equation}
\end{definition}

\hfill

\section{Variational Accelerated Optimization} \label{sec: Optimization Section}

\subsection{Riemannian Bregman Hamiltonian Approach}
\cite{Duruisseaux2021Riemannian} formulated a variational framework for the minimization of any $\lambda$-weakly-quasi-convex function $f : \mathcal{Q} \rightarrow \mathbb{R}$, via a $p$-Bregman Lagrangian  $\mathcal{L}_{p} :   T\mathcal{Q} \times \mathbb{R} \rightarrow  \mathbb{R} $ and a corresponding $p$-Bregman Hamiltonian $\mathcal{H}_{p} :   T^*\mathcal{Q} \times \mathbb{R} \rightarrow  \mathbb{R} $ for $p>0$ of the form
\begin{equation} \label{BregmanLGeneral}
	\mathcal{L}_{p}(X,V,t) = \frac{t^{\frac{\zeta}{\lambda} p +1}}{2p} \langle V , V\rangle  - Cpt^{\left(\frac{\zeta}{\lambda} +1\right)p-1} f(X),  \end{equation} 
\begin{equation} \label{BregmanHGeneral}
	\mathcal{H}_{p}(X,R,t)= \frac{p}{2t^{\frac{\zeta}{\lambda}  p +1}} \llangle R , R\rrangle + Cpt^{\left(\frac{\zeta}{\lambda}  +1\right)p-1} f(X), 
\end{equation}
where $\zeta$ is given by equation \eqref{eq: zeta}. \cite{Duruisseaux2021Riemannian} showed that solutions to the  $p$-Bregman Euler--Lagrange equations converge to a minimizer of $f$ at a convergence rate of $\mathcal{O}(1/t^p)$, under suitable assumptions. 

Furthermore, \cite{Duruisseaux2021Riemannian} proved that time-rescaling the $p$-Bregman dynamics via $\tau(t) =t^{\mathring{p}/p}$ yields the $\mathring{p}$-Bregman dynamics. Thus, the entire subfamily of Bregman trajectories indexed by the parameter $p$ can be obtained by speeding up or slowing down along the Bregman curve corresponding to any value of $p$. Inspired by the computational efficiency of the approach introduced in \cite{duruisseaux2020adaptive} on vector spaces, we can exploit the time-rescaling property of the Bregman dynamics together with a carefully chosen Poincar\'e transformation to transform the $p$-Bregman Hamiltonian into an autonomous version of the $\mathring{p}$-Bregman Hamiltonian in extended phase-space, where $\mathring{p} < p$. This allows one to integrate the higher-order $p$-Bregman dynamics while benefiting from the computational efficiency of integrating the lower-order $\mathring{p}$-Bregman dynamics. Explicitly, it was shown in \cite{Duruisseaux2021Riemannian} that the use of the time rescaling $\tau(t) = t^{\mathring{p}/p}$ within the Poincar\'e transformation framework yields the Direct approach Riemannian $p$-Bregman Hamiltonian
\begin{equation*}  \label{H-R-Direct}
	\begin{aligned}
		\bar{\mathcal{H}}_{p}(\bar{Q},\bar{R})  & = 	 \frac{p\llangle R , R\rrangle  }{2\mathfrak{Q}^{\frac{\zeta}{\lambda} p +1}} + \mathfrak{R}   + Cp\mathfrak{Q}^{\left(\frac{\zeta}{\lambda} +1\right)p-1} f(Q),
	\end{aligned}
\end{equation*}
and the Adaptive Riemannian $p\rightarrow \mathring{p}$ Bregman Hamiltonian
\begin{equation*}  \label{H-R-Adaptive}
	\begin{aligned}
		\bar{\mathcal{H}}_{p \rightarrow \mathring{p}}(\bar{Q},\bar{R})  & = 	 \frac{p^2}{2\mathring{p} \mathfrak{Q}^{\frac{\zeta}{\lambda} p +\frac{\mathring{p}}{p}}}  \llangle R , R\rrangle + \frac{p}{\mathring{p}} \mathfrak{Q}^{1-\frac{\mathring{p}}{p}}  \mathfrak{R}\\ & \qquad  \qquad \qquad  +  \frac{Cp^2}{\mathring{p}} \mathfrak{Q}^{\left(\frac{\zeta}{\lambda} +1\right)p-\frac{\mathring{p}}{p}} f(Q),
	\end{aligned}
\end{equation*}
in the extended phase space defined by $\bar{Q} = \left[\begin{smallmatrix} Q\\ \mathfrak{Q} \end{smallmatrix} \right] $ and $\bar{R} = \left[ \begin{smallmatrix} R \\ \mathfrak{R} \end{smallmatrix} \right] $ where $\mathfrak{Q}=t$  and $\mathfrak{R}$ is its conjugate momentum.

On normed vector spaces, these Riemannian Bregman Hamiltonians reduce to the Bregman Hamiltonians derived in \cite{duruisseaux2020adaptive}.
The careful computational study from \cite{duruisseaux2020adaptive}  showed that time-adaptive Hamiltonian variational discretizations, which are automatically symplectic, with adaptive time-steps informed by the time invariance of the family of $p$-Bregman Hamiltonians yielded the most robust and computationally efficient optimization algorithms, outperforming fixed-timestep symplectic discretizations, adaptive-timestep non-symplectic discretizations, and Nesterov's accelerated gradient algorithm which is neither time-adaptive nor symplectic. 

\cite{Duruisseaux2021Constrained} incorporated holonomic constraints into variational integrators to constrain the numerical solution to the Riemannian manifold, but the resulting integrators were implicit, which significantly increases their cost. 
Here, we take a different approach using the fact that the Bregman Hamiltonian in the embedding space restricts to the Riemannian Bregman Hamiltonian on the Riemannian submanifold $\mathcal{Q}$, and the projection of the Bregman Hamiltonian vector field in the embedding space onto the tangent bundle $T\mathcal{Q}$ of the Riemannian submanifold recovers the Hamiltonian vector field of the Riemannian Bregman Hamiltonian. As such, we will numerically integrate the Bregman dynamics in the embedding space and use projections to force the numerical solution to lie on $\mathcal{Q}$. If projections onto the constraint manifold $\mathcal{Q}$ can be computed exactly or approximately very efficiently, we can simply project the updated position onto $\mathcal{Q}$ after every iteration. Furthermore, if projections onto tangent spaces $T_q \mathcal{Q}$ for any point $q\in \mathcal{Q}$ are also available at a low computational cost, it might sometimes be helpful to project the update vector onto $T_q \mathcal{Q}$. Projections have been found for most Riemannian manifolds of practical interest (see \cite{Absil2008,Boumal2020}). These typically involve standard matrix factorizations which can be efficiently computed using iterative methods, and if they are expensive to compute, there are usually ways to accelerate the computations via approximations.

\subsection{Riemannian Optimization Problems}  \label{sec: List of Problems}

\subsubsection{Rayleigh Quotient Minimization on the Unit Sphere}

Eigenvectors corresponding to the largest eigenvalue of a symmetric $n \times n$ matrix $A$ maximize the Rayleigh quotient $\frac{v^\top Av}{ v^\top v }$ over $\mathbb{R}^n$. Thus, a unit eigenvector corresponding to the largest eigenvalue of the matrix $A$ is a minimizer of the function $f(v) = -  v^\top Av,$ over the unit sphere $\mathcal{Q} = \mathbb{S}^{n-1}$, which can be thought of as a Riemannian submanifold with constant positive curvature $K=1$ of $\mathbb{R}^{n}$ endowed with the Riemannian metric inherited from the Euclidean inner product $g_v(u,w) = u^\top w$. A choice of projection from $\mathbb{R}^n$ to $\mathbb{S}^{n-1}$ is the rescaling $v \mapsto \frac{v}{\|v\|_2 }$. Solving the Rayleigh quotient optimization problem efficiently is challenging when the given symmetric matrix $A$ is ill-conditioned and high-dimensional. Note that an efficient algorithm that solves the above minimization problem can also be used to find eigenvectors corresponding to the smallest eigenvalue of $A$, since the eigenvalues of $A$ are the negative of the eigenvalues of $-A$.

\subsubsection{Eigenvalue and Procrustes Problems on $\text{St}(m,n)$}

When endowed with the Riemannian metric $g_X(A,B) = \text{Trace}(A^\top B)$, the Stiefel manifold 
\begin{equation} \text{St}(m,n) = \{X\in \mathbb{R}^{n\times m} | X^\top X= I_m \} ,
\end{equation}
is a Riemannian submanifold of $\mathbb{R}^{n\times m}$. The tangent space at any $X \in \text{St}(m,n)$ is given by \begin{equation}
T_X \text{St}(m,n) = \{ Z\in \mathbb{R}^{n\times m} | X^\top Z+Z^\top X=0  \} ,\end{equation}
and the orthogonal projection $P_X$ onto $T_X \text{St}(m,n)$ is given by $	P_X Z = Z - \frac{1}{2} X(X^\top Z + Z^\top X).$ We can define a projection of any matrix $\tilde{X} \in \mathbb{R}^{n\times m}$ onto $ \text{St}(m,n)$ as the solution of 
\begin{equation*}
	\argmin_{X \in \text{St}(m,n)} \|X-\tilde{X}\|_F.
\end{equation*}
From  \cite{HaLuWa2006}, the solution of this problem is given by $X = UV^\top $ where $\tilde{X} = U\Sigma V^\top $ is the Singular Value Decomposition of $\tilde{X}$ where $\Sigma$ is a square diagonal $m \times m$ matrix. The solution $X$ of this problem can also be thought of as the first component of the polar decomposition $ \tilde{X} = X S^{1/2} $ where $X  \in \text{St}(m,n)$ and $S$ is a $m\times m$ symmetric positive-definite matrix. This solution can be written in closed form as $X = \tilde{X}(\tilde{X}^\top \tilde{X})^{-1/2}$ (and $S = \tilde{X}^\top \tilde{X}$). Thus, a first projection of any given matrix $Q \in \mathbb{R}^{n\times m}$ with Singular Value Decomposition $Q = U\Sigma V^\top$ onto $ \text{St}(m,n)$ is given by
\begin{equation*}
	Q \mapsto Q(Q^\top Q)^{-1/2} \quad \text{or equivalently} \quad Q \mapsto UV^\top. 
\end{equation*}

Another method to project a matrix $Y \in \mathbb{R}^{n\times m}$ onto $\text{St}(m,n)$ is obtained via the matrix orthogonalization $Y \mapsto \text{qf}(Y)$,  which maps the matrix $Y$ to the $Q$ factor of its QR factorization $Y =QR$ where $Q\in \text{St}(m,n)$ and $R$ is an upper triangular $n\times m$ matrix with strictly positive diagonal elements \cite{Absil2008}. 

These polar decomposition and matrix orthogonalization can also be used to construct retractions on $\text{St}(m,n) $:
\begin{equation*}
	\begin{aligned}
		\mathcal{R}_X(\xi)  = (X+\xi ) (I_m + \xi^\top \xi)^{-1/2},  
		\quad	\mathcal{R}_X(\xi)  = \text{qf} (X+\xi).
	\end{aligned}
\end{equation*}

A generalized eigenvector problem consists of finding the $m$ smallest eigenvalues of a $n\times n$ symmetric matrix $A$ and corresponding eigenvectors. This problem can be formulated as a Riemannian optimization problem on the Stiefel manifold $\text{St}(m,n)$ via the Brockett cost function
\begin{equation}
	f:\text{St}(m,n) \rightarrow \mathbb{R}, \quad   f(X) = \text{Trace}(X^\top AXN),
\end{equation}
where $N = \text{diag}(\mu_1 , \ldots , \mu_m)$ for arbitrary $0 \leq \mu_1 \leq \ldots \leq \mu_m $. The columns of a global minimizer of $f$ are eigenvectors corresponding to the $m$ smallest eigenvalues of $A$ \cite{Absil2008}.  If we define $\bar{f} : \mathbb{R}^{n\times m} \rightarrow \mathbb{R}$ via $X\mapsto \bar{f}(X) = \text{Trace}(X^\top AXN),$ then $f= \bar{f} |_{ \text{St}(m,n)} $ so
\begin{equation*} \text{grad}f(X) = P_X \text{grad}\bar{f}(X) =P_X( 2AXN). \end{equation*}

The unbalanced orthogonal Procrustes problem consists of minimizing the function
\begin{equation}
	f:\text{St}(m,n) \rightarrow \mathbb{R}, \quad         f(X) = \| AX-B \|_F^2 ,
\end{equation}
on $\text{St}(m,n) $, for given matrices $A \in \mathbb{R}^{l\times n }$ and $B \in \mathbb{R}^{l\times m }$ with $l \geq n$ and $l>m$, where $\| \cdot \|_F$ is the Frobenius norm $\| X \|_F^2 = \text{Trace}(X^\top X)=  \sum_{ij}{X_{ij}^2}$. If we define $\bar{f} : \mathbb{R}^{n\times m} \rightarrow \mathbb{R}$ via $X\mapsto \bar{f}(X) = \| AX-B \|_F^2,$ then $f= \bar{f} |_{ \text{St}(m,n)} $ so
\begin{equation*} \text{grad}f(X) = P_X \text{grad}\bar{f}(X) = P_X( 2A^\top (AX-B)). \end{equation*}
The special case where $n=m$ is the balanced orthogonal Procrustes problem. In this case, $\text{St}(m,n) = O(n)$ so $\|  AX  \|_F^2 = \|  A  \|_F^2$ and thus the minimization problem becomes the problem of maximizing $\text{ Trace}(X^\top A^\top B)$ over $X\in O(n)$. In this special case, a solution is given by $X^* = UV^\top  $ where $B^\top A = U \Sigma V^\top $ where $U$ and $V$ are orthogonal square matrices is the Singular Value Decomposition of $B^\top A$, and the solution is unique provided that the matrix $B^\top A$ is nonsingular \cite{Elden1999,Golub2013}.

\subsection{Numerical Methods}

\subsubsection{Euler--Lagrange Simple Discretization} In \cite{Duruisseaux2021Riemannian}, the $p$-Bregman Euler--Lagrange equations were rewritten as a first order system of differential equations, for which a Riemannian version of a semi-implicit Euler scheme was applied to obtain the iterates presented in Algorithm \ref{Alg: Semi-Implicit Euler}, when given a $\lambda$-weakly-quasi-convex function $f : \mathcal{Q} \rightarrow \mathbb{R}$, a retraction $\mathcal{R}$ from $T\mathcal{Q}$ to $\mathcal{Q}$, constants $C,h,p>0$, and $X_0 \in \mathcal{Q}$, $V_0 \in T_{X_0} \mathcal{Q}$.

\begin{algorithm}[!ht] 
	
	$b_k  \leftarrow 1 -  \frac{\zeta p+\lambda}{ \lambda k}, \quad c_k \leftarrow C p^2 (kh)^{p-2} $ \\
	\textbf{Version I}:	$a_k \leftarrow b_k V_k - hc_k \text{ grad}f(X_k)$ \\
	\textbf{Version II}:	$a_k \leftarrow b_k V_k - h c_k \text{ grad}f(\mathcal{R}_{X_k}(hb_k V_k))$ \;
	
	$X_{k+1} \leftarrow \mathcal{R}_{X_k}(ha_k) \quad \text{and} \quad   V_{k+1} \leftarrow \Gamma_{X_k}^{X_{k+1}} a_k$ \;
	\caption{Euler--Lagrange Discretization} \label{Alg: Semi-Implicit Euler}
\end{algorithm} 

Version I of Algorithm \ref{Alg: Semi-Implicit Euler} corresponds to the update with the traditional gradient $\nabla f(X_k)$ for the semi-implicit Euler scheme, while Version II is inspired by the reformulation of Nesterov's method from \cite{Sutskever2013} that uses a corrected gradient $\nabla f(X_k +h b_kV_k)$ instead of $\nabla f(X_k)$.

\subsubsection{Hamiltonian Taylor Variational Integrators}

\noindent Hamiltonian Taylor variational integrators (HTVIs) form a class of variational integrators described in \cite{ScShLe2017}. A discrete approximate Hamiltonian is constructed by approximating the flow map and the trajectory associated with the boundary values using a Taylor method, and approximating the integral by a quadrature rule. The variational integrator is then generated by the discrete Hamilton's equations. We will use projected versions of the HTVIs constructed in \cite{duruisseaux2020adaptive}. Given a function $f : \mathcal{Q} \rightarrow \mathbb{R}$, a projection $\mathcal{P_Q}$ onto $\mathcal{Q}$, $(q_0 , r_0)  \in T^*_{q_{0}} \mathcal{Q}$, and constants $C,h,p, \mathring{p},\mathfrak{q}_0 >0$,  the Direct and Adaptive HTVIs are obtained by iterating the updates given in Algorithm \ref{Alg: HTVI}.

\begin{algorithm}[!ht] 

	
	\textbf{Adaptive HTVI} 
	\begin{equation*} \begin{aligned}
			& r_{k+1}  = r_k -  \frac{p^2}{\mathring{p}} hC \mathfrak{q}_k^{2p-\mathring{p}/p} \text{grad}f (q_k), \\
			& q_{k+1} = \mathcal{P_Q} \left(  q_k +  \frac{p^2}{\mathring{p} } h\mathfrak{q}_k^{-p-\mathring{p}/p} r_{k+1} \right) ,  \\
			&\mathfrak{q}_{k+1}   = \mathfrak{q}_k + \frac{p}{\mathring{p} } h\mathfrak{q}_k^{1-\mathring{p}/p}.
	\end{aligned}\end{equation*}
\textbf{Direct HTVI} is the special case where $\mathring{p} = p$.
	\caption{ Direct and Adaptive HTVIs} \label{Alg: HTVI}
\end{algorithm}

\subsubsection{Riemannian Gradient Descent}

Riemannian Gradient Descent is a generalization of Gradient Descent to Riemannian manifolds which involves the Riemannian gradient and a retraction. Given a function $f : \mathcal{Q} \rightarrow \mathbb{R}$ with a retraction $\mathcal{R}$ from $T\mathcal{Q}$ to $\mathcal{Q}$, $h>0$, and $X_0 \in  \mathcal{Q}$, the Riemannian Gradient Descent iterations are given by
\begin{equation}
	X_{k+1} =  \mathcal{R}_{X_k} (-h \text{ grad} f(X_k)).
\end{equation}

\hfill 

\subsection{Numerical Results}

\subsubsection{Comparison of the Adaptive and Direct approaches}

\begin{figure}[!ht] \centering
	\includegraphics[width=0.48\textwidth]{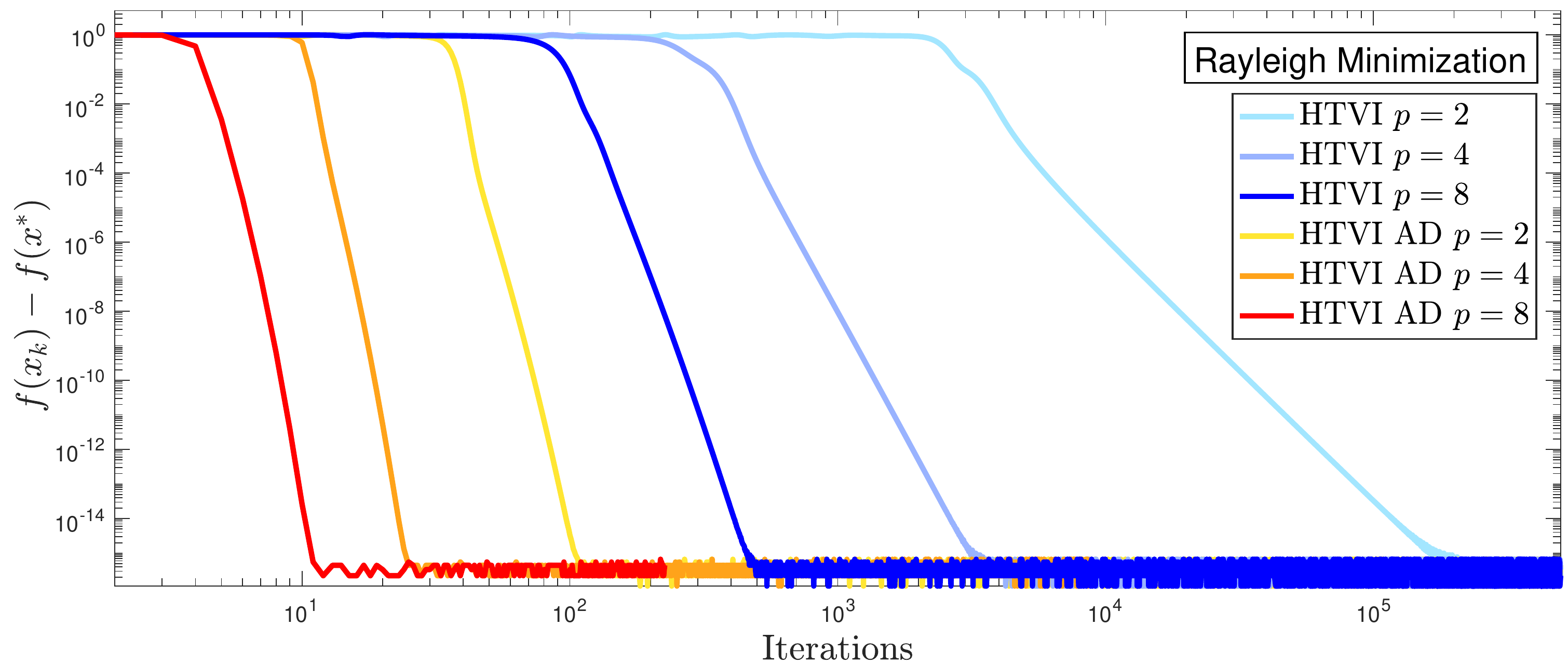}
	\caption
	{Comparison of the Direct and Adaptive (AD) projection based HTVIs with different values of the parameter $p$ and the same time-step $h = 0.01$, for the Rayleigh minimization problem on $\mathbb{S}^{n-1}$. \label{fig: HTVIpEvolution}}
\end{figure}

Numerical experiments were conducted for the Rayleigh quotient minimization problem on $\mathbb{S}^{n-1}$ with the projection based method. As was observed in \cite{duruisseaux2020adaptive} on vector spaces, Figure \ref{fig: HTVIpEvolution} shows that the Adaptive approach can be significantly more efficient than the Direct approach, and that both methods enjoy faster convergence as $p$ increases.

\subsubsection{Rayleigh minimization problem on the unit sphere  $\mathbb{S}^{n-1}$}

\begin{figure}[!ht] \centering 
	\includegraphics[width=0.48\textwidth]{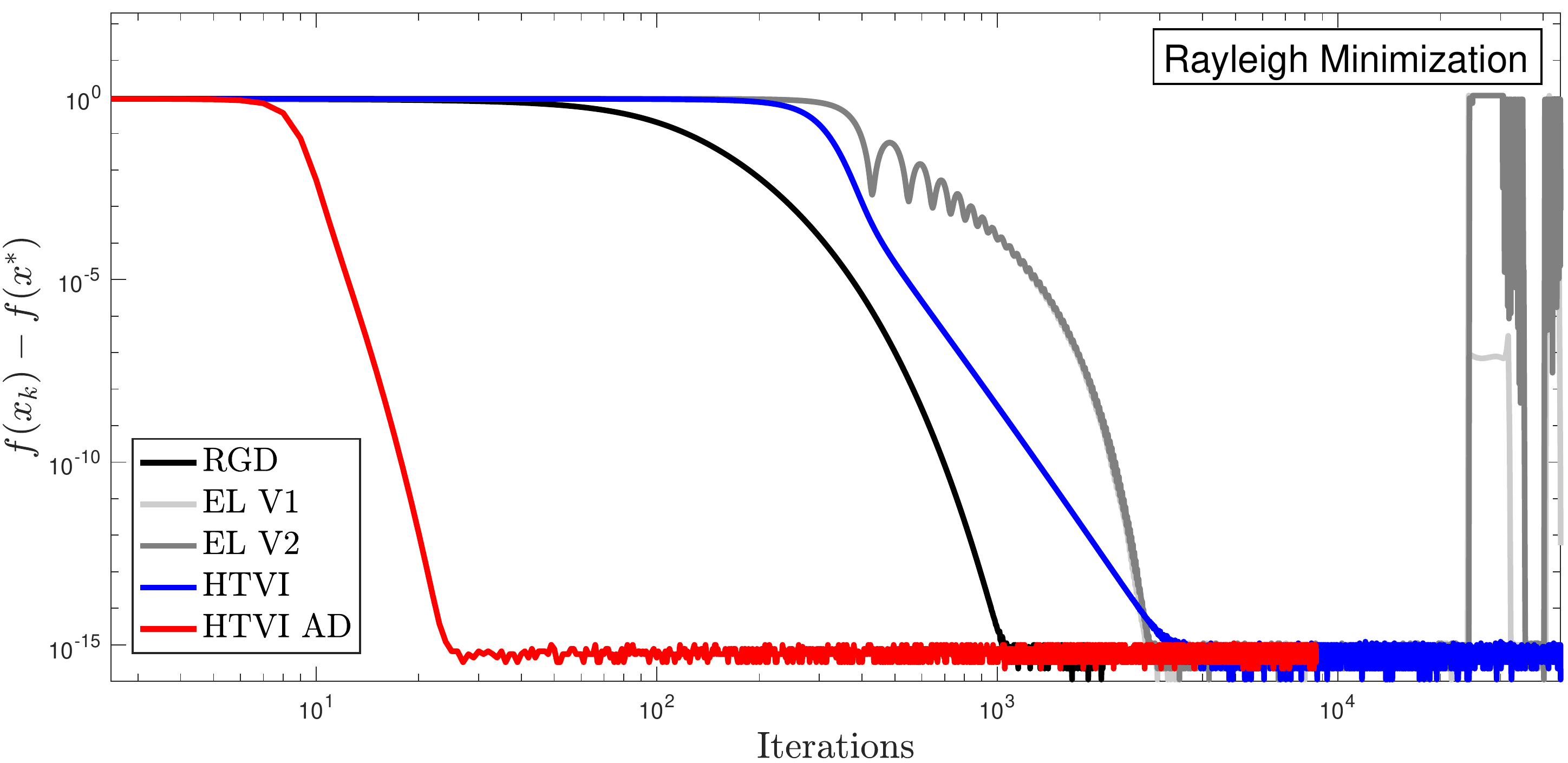}
	\caption{Comparison of the Direct and Adaptive (AD) projection based HTVIs with the Riemannian Gradient Descent (RGD) method and the Euler--Lagrange discretizations (EL V1 and EL V2), with $p=4$  and the same time-step $h$.  \label{fig: ComparisonAll} }
\end{figure}

It was noted in \cite{Duruisseaux2021Riemannian} that although higher values of $p$ in Algorithm \ref{Alg: Semi-Implicit Euler} result in provably faster rates of convergence, they also appear to be more prone to stability issues under numerical discretization, which can cause the numerical optimization algorithm to diverge. Numerical experiments in \cite{duruisseaux2020adaptive} showed that on normed vector spaces, geometric discretizations which respect the time-rescaling invariance and symplecticity of the Bregman Hamiltonian flows were substantially less prone to these stability issues, and were therefore more robust, reliable, and computationally efficient. This was one of the motivations to develop time-adaptive variational integrators for the Bregman Hamiltonians.  Numerical experiments were conducted for the Rayleigh quotient minimization problem on $\mathbb{S}^{n-1}$, and the results are presented in Figure~\ref{fig: ComparisonAll}.

The Adaptive HTVI clearly outperforms the other algorithms for this problem. Note that the Euler--Lagrange discretizations suffer from stability issues leading to a loss of convergence (after $\approx 10^4$ iterations), due to the polynomially growing unbounded term $C p^2 (kh)^{p-2} $ paired with the $\text{grad} f$ term  to 0 which can only achieve a finite order of accuracy due to numerical roundoff error. This issue can be fixed by adding a suitable upper bound to the coefficient $C p^2 (kh)^{p-2} $ in the updates, or by stopping the iterating process once a desired convergence criterion is achieved.

\subsubsection{Optimization Problems on the Stiefel manifold $\text{St}(m,n)$}

\begin{figure}[!ht]
	\centering
	\begin{minipage}[b]{0.48\textwidth}
		\includegraphics[width=\textwidth]{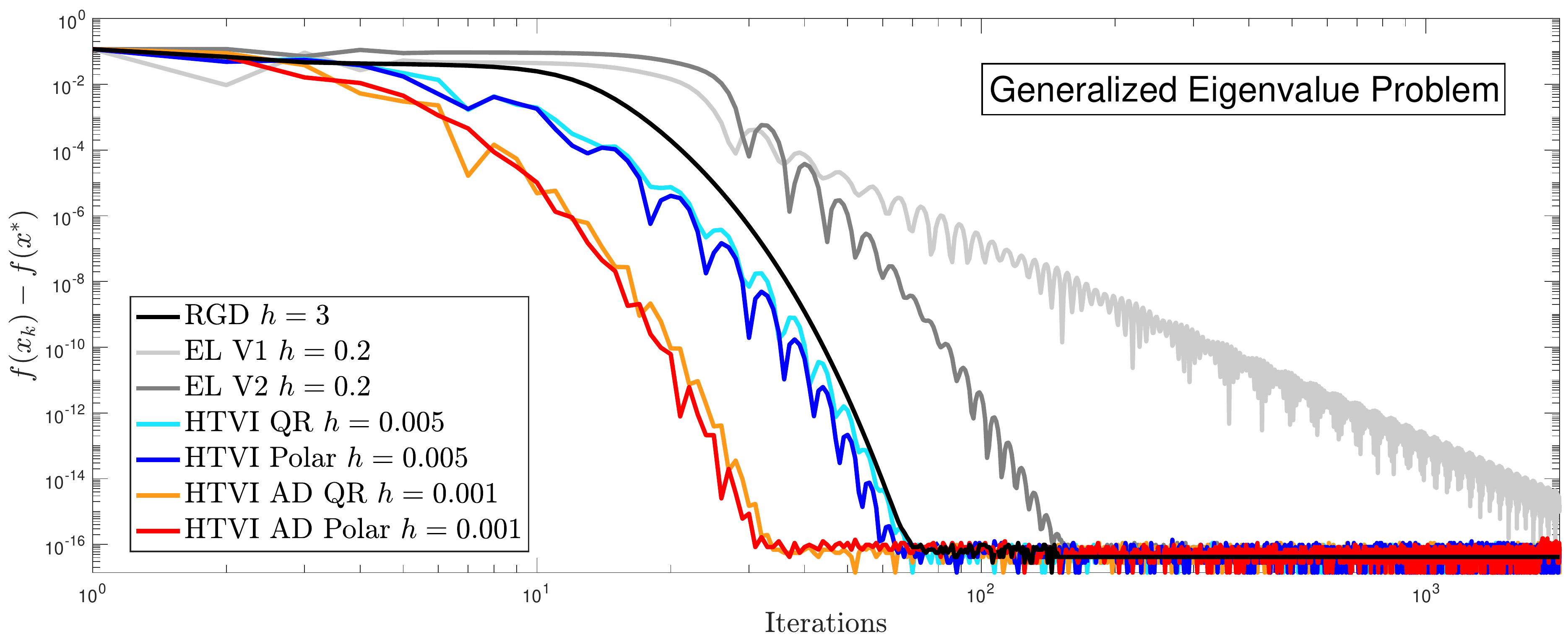}
	\end{minipage}
	\\ 
	\begin{minipage}[b]{0.48\textwidth}
		\includegraphics[width=\textwidth]{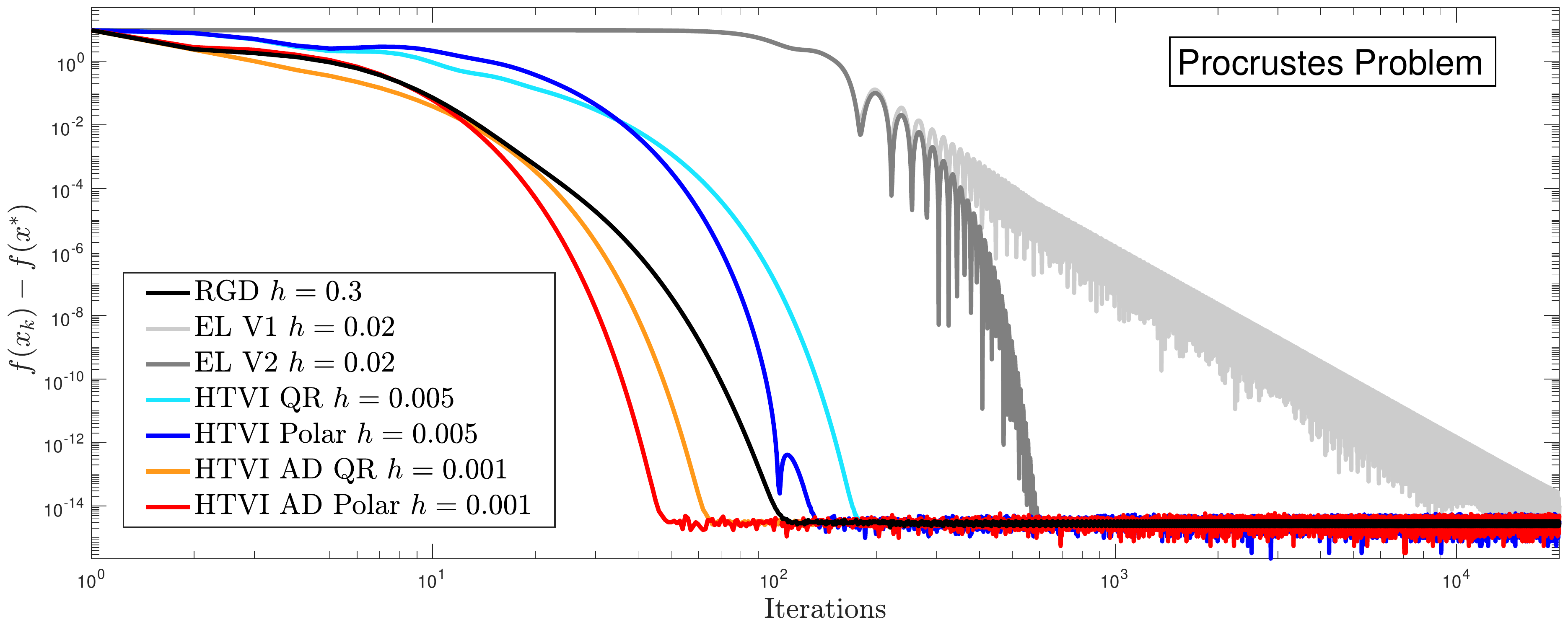}
	\end{minipage}
	\\ 
	\begin{minipage}[b]{0.48\textwidth}
		\includegraphics[width=\textwidth]{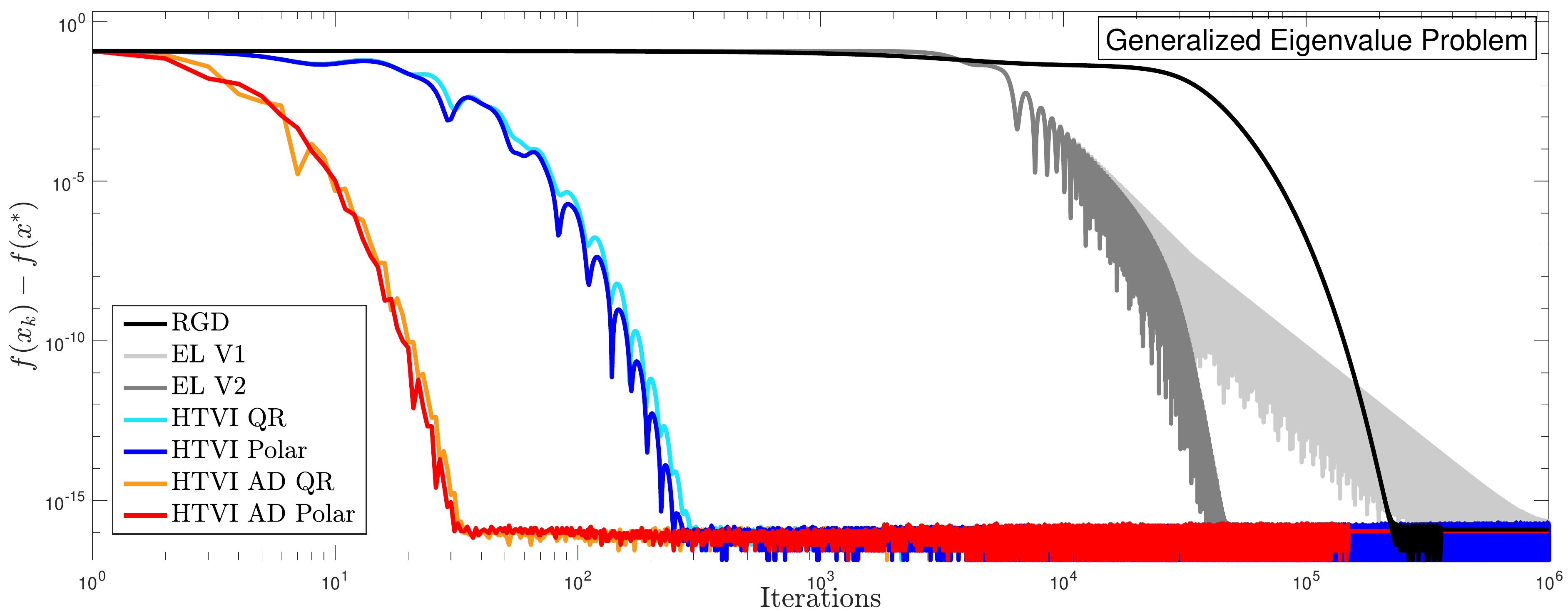}
	\end{minipage} 
	\\ 
	\begin{minipage}[b]{0.48\textwidth}
		\includegraphics[width=\textwidth]{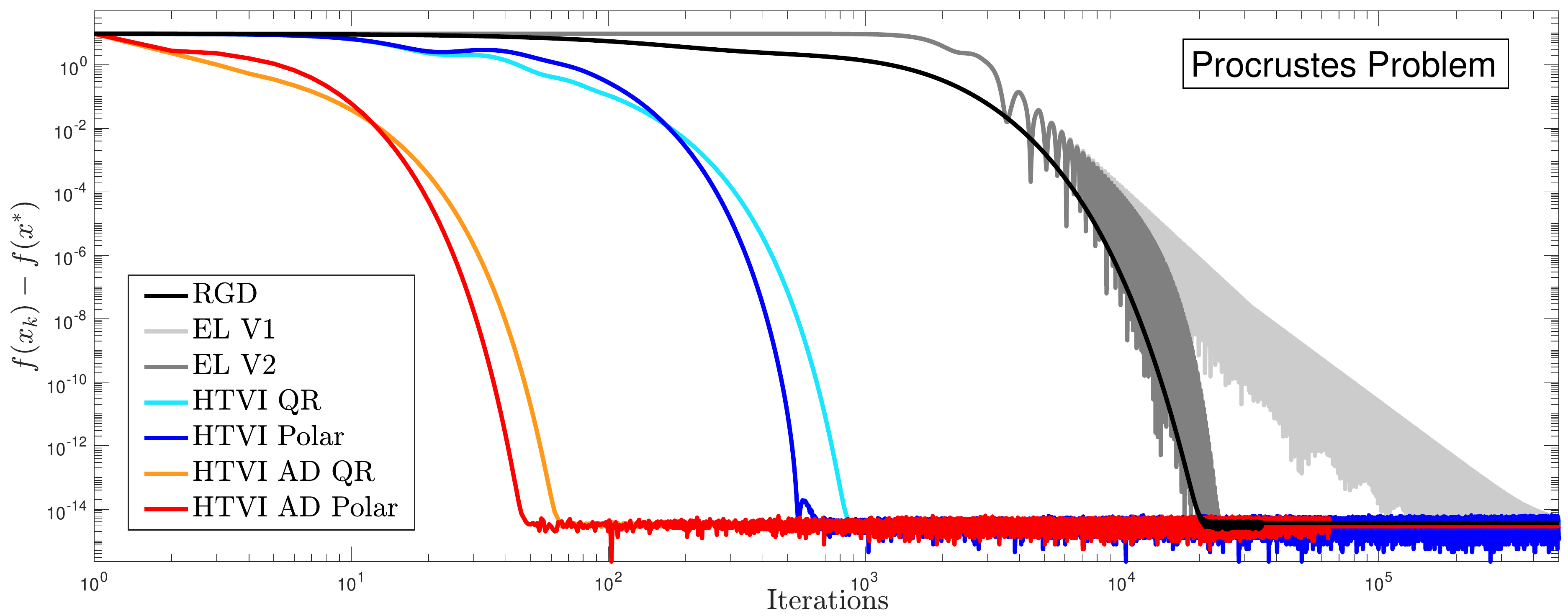}
	\end{minipage} 
	\caption{Comparison of the Direct and Adaptive (AD) Type II HTVIs with the Riemannian Gradient Descent (RGD) method and the Euler--Lagrange discretizations (EL V1 and EL V2) with $p=5$ with different time steps (top two plots) and with the same time-step $h = 0.001$ (bottom two plots), for the generalized eigenvalue and Procrustes problems on $\text{St}(m,n)$.  \label{fig: Stiefel} }
\end{figure}

Numerical experiments were conducted for the generalized eigenvalue and Procustes problems on $\text{St}(m,n)$ to observe how the projection based HTVIs compare to the Euler--Lagrange discretizations from \cite{Duruisseaux2021Riemannian} and the standard Riemannian gradient descent. The results are presented in Figure \ref{fig: Stiefel}. Note that for certain instances of the Procrustes problem with certain initial values $X_0 \in \text{St}(m,n)$, all the algorithms converged to a local minimizer which was not a global minimizer.

The projection based Adaptive HTVIs clearly outperform their Direct approach counterparts, Riemannian gradient descent and both versions of the Euler--Lagrange discretization in terms of number of iterations required, when all the algorithms are implemented with the same time step (see the two bottom plots in Figure \ref{fig: Stiefel}). As can be seen from the top two plots in Figure \ref{fig: Stiefel}, the Adaptive HTVIs are still the best performing algorithms, even when larger time steps are taken for the other algorithms and in particular even when the Riemannian gradient descent algorithm has been tuned optimally. Note that both the Euler--Lagrange discretizations and Hamiltonian variational integrators suffered from the numerical roundoff issue described in the previous subsection, but this issue was resolved by adding a suitable upper bound to the ever-growing problematic coefficient in the updates.

Our numerical experiments do not suggest that there is a clear benefit in using the polar decomposition based projection over the matrix orthogonalization, or vice versa. Both projection strategies led to very efficient algorithms for Riemannian accelerated optimization with seemingly similar performance and stability properties. Computing the QR decomposition of a $n \times m$ matrix via the standard Householder QR algorithm requires approximately $2m^2( n -  m/3)$ floating point operations, while computing the Singular Value Decomposition of a $n \times m$ is more expensive and often relies on intermediate QR decompositions \cite{Trefethen1997}. Thus, these operations can become very costly as the dimension of the problem becomes large, in which case it might be beneficial to use approximate QR decompositions and Singular Value Decompositions. For instance, the projection based on the polar decomposition $Q \mapsto Q(Q^\top Q)^{-1/2}$ can be rewritten as $	Q \mapsto Q(I_m + ( Q^\top Q - I_m))^{-1/2}$, and provided the distance away from the Stiefel manifold is sufficiently small, the norm of $E = (Q^\top Q - I_m)$ is small and we can approximate the projection by truncating its series expansion
\begin{equation*}
	Q(I_m + D)^{-1/2} = Q \left(I_m - \frac{1}{2}D + \frac{3}{8} D^2- \frac{5}{16}D^3 + \ldots \right).
\end{equation*}

We also tested the projection algorithm against the implicit algorithm from \cite{Duruisseaux2021Constrained} on the same optimization problems on $\mathbb{S}^{n-1}$ and $\text{St}(m,n)$ Although both algorithms produced very similar graphs for the error as a function of the iteration number, the explicit nature of our projection algorithm made every iteration significantly faster and overall the running time was reduced by several orders of magnitude, even on low-dimensional problems (for instance, 3 orders of magnitude on $\mathbb{S}^{5-1}$ and $\text{St}(3,2)$, and 4 orders of magnitude on $\mathbb{S}^{100-1}$). Note that the projection algorithm was also easier to implement and tune than the implicit algorithm.

\section{Conclusion}

Motivated by the observation made in the normed space setting in \cite{duruisseaux2020adaptive} that a careful use of adaptivity and symplecticity within the variational formulation of accelerated optimization could result in a significant gain in computational efficiency, discrete variational integrators incorporating holonomic constraints were constructed in \cite{Duruisseaux2021Constrained} within the variational framework for Riemannian accelerated optimization of \cite{Duruisseaux2021Riemannian}. The resulting algorithms performed well in terms of number of iterations required to achieve convergence but were implicit which can lead to high computational costs. In this paper, we saw that the gain in computational efficiency is preserved when the constraints are enforced via projections instead of being incorporated directly into the variational principles, and that the explicit nature of the resulting algorithms makes every iteration significantly faster and easier to tune than for the implicit algorithms from \cite{Duruisseaux2021Constrained}. As a consequence, if projections onto the constraint manifold can be computed efficiently, these projection based variational integrators form a class of efficient explicit algorithms for Riemannian accelerated optimization, and we believe that these algorithms are the most efficient methods to date which exploit the variational framework from \cite{Duruisseaux2021Riemannian}.

Although the Whitney and Nash Embedding Theorems imply that there is no loss of generality when studying Riemannian manifolds only as submanifolds of Euclidean spaces, designing intrinsic methods that would exploit and preserve the symmetries and geometric properties of the Riemannian manifold and of the problem at hand could have advantages both in terms of computation and in terms of improving our understanding of the acceleration phenomenon on Riemannian manifolds. Developing an intrinsic extension of Hamiltonian variational integrators to manifolds would require some additional work, since the current approach involves Type II/III generating functions $H_d^+(q_k, p_{k+1})$, $H_d^-(p_k, q_{k+1})$, which depend on the position at one boundary point, and the momentum at the other boundary point. This does not make intrinsic sense on a manifold, since one needs the base point in order to specify the corresponding cotangent space, and instead one should ideally consider a construction based on the more general discrete Dirac mechanics~\cite{LeOh2008}.

It would be desirable to have some convergence guarantees, but proving that the discrete time algorithms perform analogously to the continuous dynamics is far from direct, as the $\mathcal{O}(1/t^p)$ convergence rate for the continuous-time dynamics conflicts with the $\mathcal{O}(1/k^2)$ Nesterov barrier theorem for discrete-time algorithms. Some shadowing results can be obtained for certain Riemannian optimization algorithms when the associated dynamical system is uniformly contracting. However, momentum methods such as the ones presented here, are notoriously non-descending and heavily oscillatory, and lack contraction as a result. Even on vector spaces, obtaining theoretical guarantees was a challenging task, achieved in \cite{Jadbabaie2018} under additional assumptions. Generalizing these results to Riemannian manifolds would be much more challenging than a trivial generalization of these results since vector space operations and objects have to be replaced by their Riemannian generalizations which involve geodesics, parallel transport, covariant derivatives, Riemannian exponentials and logarithms.



\newpage

%
%


\bibliographystyle{named}
\small
\bibliography{ExplicitBib}

\def\cprime{$'$}
\begin{thebibliography}{}

\bibitem[\protect\citeauthoryear{Absil \bgroup \em et al.\egroup
  }{2008}]{Absil2008}
P.~A. Absil, R.~Mahony, and R.~Sepulchre.
\newblock {\em Optimization Algorithms on Matrix Manifolds}.
\newblock Princeton University Press, Princeton, NJ, 2008.

\bibitem[\protect\citeauthoryear{Ahn and Sra}{2020}]{Sra2020}
K.~Ahn and S.~Sra.
\newblock From {N}esterov's estimate sequence to {R}iemannian acceleration.
\newblock In {\em Proceedings of Thirty Third Conference on Learning Theory},
  volume 125 of {\em Proceedings of Machine Learning Research}, pages 84--118.
  PMLR, 09--12 Jul 2020.

\bibitem[\protect\citeauthoryear{Alimisis \bgroup \em et al.\egroup
  }{2020}]{Alimisis2020-1}
F.~Alimisis, A.~Orvieto, G.~B\'ecigneul, and A.~Lucchi.
\newblock Practical accelerated optimization on {R}iemannian manifolds, 2020.

\bibitem[\protect\citeauthoryear{Benettin and Giorgilli}{1994}]{Benettin1994}
G.~Benettin and A.~Giorgilli.
\newblock On the {H}amiltonian interpolation of near-to-the identity symplectic
  mappings with application to symplectic integration algorithms.
\newblock {\em Journal of Statistical Physics}, 74:1117--1143, 03 1994.

\bibitem[\protect\citeauthoryear{Boumal}{2020}]{Boumal2020}
N.~Boumal.
\newblock An introduction to optimization on smooth manifolds, 2020.
\newblock Available online at http://www.nicolasboumal.net/book.

\bibitem[\protect\citeauthoryear{Duruisseaux and
  Leok}{2021a}]{Duruisseaux2021Constrained}
V.~Duruisseaux and M.~Leok.
\newblock Accelerated optimization on {R}iemannian manifolds via discrete
  constrained variational integrators.
\newblock 2021.

\bibitem[\protect\citeauthoryear{Duruisseaux and
  Leok}{2021b}]{Duruisseaux2021Riemannian}
V.~Duruisseaux and M.~Leok.
\newblock A variational formulation of accelerated optimization on {R}iemannian
  manifolds.
\newblock 2021.

\bibitem[\protect\citeauthoryear{Duruisseaux \bgroup \em et al.\egroup
  }{2021}]{duruisseaux2020adaptive}
V.~Duruisseaux, J.~Schmitt, and M.~Leok.
\newblock Adaptive {H}amiltonian variational integrators and applications to
  symplectic accelerated optimization.
\newblock {\em SIAM Journal on Scientific Computing}, 43(4):A2949--A2980, 2021.

\bibitem[\protect\citeauthoryear{Eld{\'{e}}n and Park}{1999}]{Elden1999}
L.~Eld{\'{e}}n and H.~Park.
\newblock A {P}rocrustes problem on the {S}tiefel manifold.
\newblock {\em Numerische Mathematik}, 82(4):599--619, 1999.

\bibitem[\protect\citeauthoryear{Golub and Van~Loan}{2013}]{Golub2013}
G.~H. Golub and C.~F. Van~Loan.
\newblock {\em Matrix Computations}.
\newblock Johns Hopkins Studies in the Mathematical Sciences. Johns Hopkins
  University Press, 2013.

\bibitem[\protect\citeauthoryear{Hairer \bgroup \em et al.\egroup
  }{2006}]{HaLuWa2006}
E.~Hairer, C.~Lubich, and G.~Wanner.
\newblock {\em Geometric {N}umerical {I}ntegration}, volume~31 of {\em Springer
  Series in Computational Mathematics}.
\newblock Springer-Verlag, Berlin, second edition, 2006.

\bibitem[\protect\citeauthoryear{Lall and West}{2006}]{LaWe2006}
S.~Lall and M.~West.
\newblock Discrete variational {H}amiltonian mechanics.
\newblock {\em J. Phys. A}, 39(19):5509--5519, 2006.

\bibitem[\protect\citeauthoryear{Lee}{2018}]{Lee2019}
{J. M.} Lee.
\newblock {\em Introduction to {R}iemannian Manifolds}, volume 176 of {\em
  Graduate Texts in Mathematics}.
\newblock Springer, Cham, second edition, 2018.

\bibitem[\protect\citeauthoryear{Leok and Ohsawa}{2011}]{LeOh2008}
M.~Leok and T.~Ohsawa.
\newblock Variational and geometric structures of discrete {D}irac mechanics.
\newblock {\em Found. Comput. Math.}, 11(5):529--562, 2011.

\bibitem[\protect\citeauthoryear{Leok and Zhang}{2011}]{LeZh2011}
M.~Leok and J.~Zhang.
\newblock Discrete {H}amiltonian variational integrators.
\newblock {\em IMA Journal of Numerical Analysis}, 31(4):1497--1532, 2011.

\bibitem[\protect\citeauthoryear{Liu \bgroup \em et al.\egroup
  }{2017}]{Liu2017}
Y.~Liu, F.~Shang, J.~Cheng, H.~Cheng, and L.~Jiao.
\newblock Accelerated first-order methods for geodesically convex optimization
  on {R}iemannian manifolds.
\newblock In {\em NeurIPS}, volume~30, pages 4868--4877, 2017.

\bibitem[\protect\citeauthoryear{Marsden and West}{2001}]{MaWe2001}
J.~E. Marsden and M.~West.
\newblock Discrete mechanics and variational integrators.
\newblock {\em Acta Numer.}, 10:357--514, 2001.

\bibitem[\protect\citeauthoryear{Nash}{1956}]{Nash1956}
J.~Nash.
\newblock The imbedding problem for {R}iemannian manifolds.
\newblock {\em Annals of Mathematics}, 63(1):20--63, 1956.

\bibitem[\protect\citeauthoryear{Nesterov}{1983}]{Nes83}
Y.~Nesterov.
\newblock A method of solving a convex programming problem with convergence
  rate $\mathcal{O}(1/k^2)$.
\newblock {\em Soviet Mathematics Doklady}, 27(2):372--376, 1983.

\bibitem[\protect\citeauthoryear{Reich}{1999}]{Re1999}
S.~Reich.
\newblock Backward error analysis for numerical integrators.
\newblock {\em SIAM J. Numer. Anal.}, 36:1549--1570, 1999.

\bibitem[\protect\citeauthoryear{Schmitt and Leok}{2017}]{ScLe2017}
J.~M. Schmitt and M.~Leok.
\newblock Properties of {H}amiltonian variational integrators.
\newblock {\em IMA Journal of Numerical Analysis}, 38(1):377--398, 03 2017.

\bibitem[\protect\citeauthoryear{Schmitt \bgroup \em et al.\egroup
  }{2018}]{ScShLe2017}
J.~M. Schmitt, T.~Shingel, and M.~Leok.
\newblock {L}agrangian and {H}amiltonian {T}aylor variational integrators.
\newblock {\em {BIT} Numerical Mathematics}, 58:457--488, 2018.

\bibitem[\protect\citeauthoryear{Su \bgroup \em et al.\egroup
  }{2016}]{SuBoCa16}
W.~Su, S.~Boyd, and E.~Candes.
\newblock A differential equation for modeling {N}esterov's {A}ccelerated
  {G}radient method: theory and insights.
\newblock {\em Journal of Machine Learning Research}, 17(153):1--43, 2016.

\bibitem[\protect\citeauthoryear{Sutskever \bgroup \em et al.\egroup
  }{2013}]{Sutskever2013}
I.~Sutskever, J.~Martens, G.~Dahl, and G.~Hinton.
\newblock On the importance of initialization and momentum in deep learning.
\newblock In {\em Proceedings of the 30th International Conference on
  International Conference on Machine Learning - Volume 28}, ICML'13, pages
  1139--1147, Atlanta, GA, USA, 2013.

\bibitem[\protect\citeauthoryear{Trefethen and Bau}{1997}]{Trefethen1997}
L.N. Trefethen and D.~Bau.
\newblock {\em Numerical Linear Algebra}.
\newblock Other Titles in Applied Mathematics. SIAM, 1997.

\bibitem[\protect\citeauthoryear{Whitney}{1944a}]{Whitney1944_2}
H.~Whitney.
\newblock The self-intersections of a smooth $n$-manifold in $2n$-space.
\newblock {\em Annals of Mathematics}, 45(2):220--246, 1944.

\bibitem[\protect\citeauthoryear{Whitney}{1944b}]{Whitney1944_1}
H.~Whitney.
\newblock The singularities of a smooth $n$-manifold in $(2n - 1)$-space.
\newblock {\em Annals of Mathematics}, 45(2):247--293, 1944.

\bibitem[\protect\citeauthoryear{Wibisono \bgroup \em et al.\egroup
  }{2016}]{WiWiJo16}
A.~Wibisono, A.~Wilson, and M.~Jordan.
\newblock A variational perspective on accelerated methods in optimization.
\newblock {\em Proceedings of the National Academy of Sciences},
  113(47):E7351--E7358, 2016.

\bibitem[\protect\citeauthoryear{Zhang and Sra}{2016}]{Sra2016}
H.~Zhang and S.~Sra.
\newblock First-order methods for geodesically convex optimization.
\newblock In {\em 29th Annual Conference on Learning Theory}, pages 1617--1638,
  2016.

\bibitem[\protect\citeauthoryear{Zhang and Sra}{2018}]{Sra2018}
H.~Zhang and S.~Sra.
\newblock An estimate sequence for geodesically convex optimization.
\newblock In {\em Proceedings of the 31st Conference On Learning Theory},
  volume~75 of {\em Proceedings of Machine Learning Research}, pages
  1703--1723, Jul 2018.

\bibitem[\protect\citeauthoryear{Zhang \bgroup \em et al.\egroup
  }{2018}]{Jadbabaie2018}
J.~Zhang, A.~Mokhtari, S.~Sra, and A.~Jadbabaie.
\newblock Direct {R}unge-{K}utta discretization achieves acceleration.
\newblock In {\em Advances in Neural Information Processing Systems},
  volume~31. Curran Associates, Inc., 2018.

\end{thebibliography}

\end{document}